\documentclass[twoside]{amsart}
\copyrightinfo{}{}

\usepackage{amssymb}
\usepackage{amsmath}
\usepackage{cite}
\usepackage{mathrsfs}
\usepackage{booktabs}

\makeatletter

\def\myarrow#1{\,\ooalign{\raise.8ex\hbox{\ \small#1}\crcr%
           $\xrightarrow{\phantom{\hbox{\ \small#1}}}$}\,}
\newcommand{\lra}{\begin{picture}(80,17)
                  \put(10,10){\vector(1,0){60}}
                  \end{picture}}
\newcommand{\lraa}[1]{\begin{picture}(80,50)      
                      \put(  10,  10){\vector(1,0){60}}
                      \put(  10,  27){\makebox[6mm]{$\scriptstyle #1$}}
                      \end{picture} }
\newcommand{\lramaps}{\begin{picture}(80,17)
                  \put(10,10){\vector(1,0){60}}
                  \put(10,0){\line(0,1){20}}
                  \end{picture}}
\newcommand{\lraiso}{\begin{picture}(80,30)
                     \put(   0,   0){$\lra$}
                     \put(  25,  15){$\scriptstyle\sim$}
                     \end{picture} }
\newcommand{\lraisoa}[1]{\begin{picture}(80,40)
                         \put(   0,   0){$\lra$}
                         \put(  10,  25){\makebox[6mm]{$\scriptstyle #1$}}
                         \put(  25,  -8){$\scriptstyle\sim$}
                         \end{picture} }

\let\atop\@@atop

\def\And{\text{\ and\ }}

\def\For{\text{\ for\ }}
\def\ForAll{\text{\ for all\ }}

\def\If{\text{\ if\ }}

\def\Or{\text{\ or\ }}
\def\Otherwise{\text{\ otherwise}}

\def\th{{\text{th}}}

\def\Set[#1]#2|#3|{\Big\{\ #2\ \Big| \
            \vcenter{\hsize #1mm\centering#3}\Big\}}

{\catcode`\|=\active
  \gdef\set#1{\mathinner{\lbrace\,{\mathcode`\|"8000%
                                   \let|\midvert #1}\,\rbrace}}
}
\def\midvert{\egroup\mid\bgroup}

\def\Number#1{\refstepcounter{equation}
              \leqno(\theequation)\if*#1%
              \else\def\@currentlabel{{\rm\theequation}}\label{#1}%
              \fi}
\def\Dag{\ifmmode\leqno(\dag)\else$(\dag\)$\fi}
\def\DDag{\ifmmode\leqno(\ddag)\else$(\ddag\)$\fi}

\numberwithin{equation}{section}

\swapnumbers

\newtheorem{Proposition}[equation]{Proposition}
\newtheorem{Lemma}[equation]{Lemma}
\newtheorem{Corollary}[equation]{Corollary}

\theoremstyle{remark}

\newtheorem{Example}[equation]{Example}

\newenvironment{Point}[2]%
  {\ifx*#2\let\pointlabel\relax\else\def\pointlabel{#2}\fi
   \refstepcounter{equation}\trivlist
   \item[\hskip\labelsep\theequation\bf
         \ifx\pointlabel\relax\else\space\pointlabel\space\fi]
   \ignorespaces #1
  }{\relax}

\def\Sum{\displaystyle\sum}

\def\){\big)}
\def\({\big(}
\let\>\rangle
\let\<\langle

\let\ss\subseteq
\let\0\varnothing              

\let\realb@r\bar
\let\bar\overline

\let\gedom\trianglerighteq

\def\diag{\mathop{\rm diag}\nolimits}

\def\Ind{\mathop{\rm Ind}\nolimits}

\def\GL{\mathop{\rm GL}\nolimits}
\def\Kern{\mathop{\rm Kern}\nolimits}

\def\F{{\mathbb F}}
\def\N{{\mathbb N}}
\def\Z{{\mathbb Z}}
\def\Q{{\mathbb Q}}

\def\map#1#2{\,{:}\,#1\!\longrightarrow\!#2}
\def\mapsto{\!\longmapsto\!}

\keywords{Specht modules, elementary divisors, hook partitions.}

\title{Elementary divisors of Specht modules}

\author{Matthias K\"unzer}
\address{\begin{tabular}{l@{}}
         Abteilung Reine Mathematik,\\
         Universit\"at Ulm, D-89069 Ulm.\phantom{xxxxx}
                   \hfill\qquad
                   {\tt kuenzer@mathematik.uni-ulm.de} \\
         Germany.\hfill\qquad
                   {\tt www.mathematik.uni-ulm.de/ReineM/kuenzer}
         \end{tabular}
}
\author{Andrew Mathas}
\address{\begin{tabular}{l@{}}
            School of Mathematics and Statistics F07\\  
            University of Sydney,
            Sydney N.S.W. 2006.\hfill\qquad
                        {\tt mathas@maths.usyd.edu.au}\\  
            Australia.\hfill\qquad
                    {\tt www.maths.usyd.edu.au/u/mathas/}
         \end{tabular}
}

\usepackage[enableskew,vcentermath]{youngtab}
\def\tab(#1|#2){\Big(\,\young(#1)\,\Big|\,\young(#2)\,\Big)}
\def\widetab(#1){\vcenter{\tiny%
   \offinterlineskip\halign{\vrule##&&\ \hfil$##$\hfil\ &\vrule##\cr%
         \noalign{\hrule} height 5.5pt depth 1.5pt &#1&\cr
         \noalign{\hrule}
}}}

\def\len(#1){\ell(#1)}

\def\D{\mathcal D}
\def\H{\mathscr H}

\def\Std{\operatorname{\rm Std}}
\def\Sym{\mathfrak S}
\def\ZZ{{\Z[q,q^{-1}]}}

\def\a{\mathfrak a}
\def\b{\mathfrak b}
\def\s{\mathfrak s}
\def\t{\mathfrak t}
\begin{document}
\maketitle

\section{Introduction and statement of main results}

The irreducible representations of the symmetric groups and their
Iwahori-Hecke algebras have been classified and constructed by
James~\cite{James:SnIrr} and Dipper and James~\cite{DJ:reps}, yet
simple properties of these modules, such as their dimensions, are
still not known. Every irreducible representation of these algebras is
constructed by quotienting out the radical of a bilinear form on a
particular type of module, known as a Specht module. The bilinear
forms on the Specht modules are the objects of our study.  

One way of determining the dimension of the simple modules would be to
first find the elementary divisors of its Gram matrix over
$\Z[q,q^{-1}]$ and then specialize. This would also give the
dimensions of the subquotients of the Jantzen filtrations of the
Specht modules over an arbitrary field; see \cite{JM:Jantzen}.  In
general, such an approach is not possible because, as Andersen has
shown, Gram matrices need not be diagonalizable over $\Z[q,q^{-1}]$;
see~\cite[Remark~5.11]{Andersen:Tilting}. We also give some examples
of non--diagonalizable Specht modules in section~\ref{SecCounterEx}.

Let $G(\lambda)$ be the Gram matrix of the Specht module $S(\lambda)$.
Then the first result in this paper shows that $G(\lambda)$ is
diagonalizable if and only if $G(\lambda')$ is diagonalizable,
where~$\lambda'$ is the partition conjugate to $\lambda$. Moreover, if
$G(\lambda)$ is divisibly diagonalizable (that is, $G(\lambda)$ is
equivalent to a diagonal matrix $\diag(d_1,\dots,d_m)$ such that $d_i$
divides~$d_{i+1}$, for $1\le i<m$), then so is $G(\lambda')$. In this
case we can speak of elementary divisors and we show how the
elementary divisors of $G(\lambda)$ and $G(\lambda')$ determine each
other. This is a $q$--analogue of the corresponding result for the
symmetric group~\cite{KN}.

We next consider the elementary divisors for the hook partitions. We
show that when $\lambda = (n-k,1^k)$, for $0\le k<n$, the Gram matrix
$G(\lambda)$ is always divisibly diagonalizable over $\Z[q,q^{-1}]$,
and we determine the elementary divisors. Again, this is a
$q$--analogue of the corresponding result for the symmetric
groups~\cite{KN}, however, the proof in the Hecke algebra case is more
involved and requires some interesting combinatorics.

\section{The Hecke algebra and permutation modules}

Fix a positive integer $n$ and let $\Sym_n$ be the symmetric group of
degree $n$. 

Let $R$ be a commutative domain and let $q$ be an
invertible element in $R$. 

The Iwahori--Hecke algebra of $\Sym_n$ with
parameter $q$ is the unital associative algebra $\H$ with generators
$T_1,T_2,\dots,T_{n-1}$ and relations
$$
\begin{array}{rcll}
  (T_i-q)(T_i+1)  & = & 0                   & \For 1\le i < n,\\
  T_i T_j         & = & T_j T_i             & \For 1\le i < j-1 < n-1,\\
  T_i T_{i+1} T_i & = & T_{i+1} T_i T_{i+1} & \For 1\le i < n-1 \; .
\end{array}
$$

Let $r_i=(i,i+1)$, for $i=1,2,\dots,n-1$. Then
$\{r_1,r_2,\dots,r_{n-1}\}$ generate $\Sym_n$ (as a Coxeter group).
If $w\in\Sym_n$ then $w = r_{i_1}\cdots r_{i_k}$ for some $i_j$ with $1\le i_j < n$.
The word $w = r_{i_1}\dots r_{i_k}$ is {\sf reduced} if $k$ is minimal;
in this case we say that $w$ has {\sf length} $k$ and we define $\len(w)=k$. 

If $r_{i_1}\dots r_{i_k}$ is reduced then we set $T_w = T_{i_1}\cdots T_{i_k}$. 
Then $T_w$ is independent of the choice of reduced
expression for $w$; see, for example, \cite[1.11]{M:IHA}. Furthermore, 
$\H$ is free as an $R$--module with basis $\set{T_w|w\in\Sym_n}$.

A {\sf composition} $\mu$ of $n$ is a sequence of non--negative
integers $(\mu_1,\mu_2,\dots)$ that sum to $n$. If, in addition,
$\mu_1\ge\mu_2\geq\dots$, then $\mu$ is a {\sf partition} of $n$.

Let $\mu$ be a composition of $n$ and let $\Sym_\mu$ be the associated
Young subgroup. Then $\H(\Sym_\mu)=\<\,T_w\mid w\in\Sym_\mu\,\>$ is a
subalgebra of $\H$. Given a (right) $\H(\Sym_\mu)$--module~$V$, we
define the induced $\H$--module
$$
\Ind_{\H(\Sym_\mu)}^\H(V) \; =\; V\otimes_{\H(\Sym_\mu)}\H\; .
$$
Let $\D_\mu=\set{d\in\Sym_n|\len(dr_i)>\len(d)\ForAll r_i\in\Sym_\mu}$ 
be the set of distinguished right coset
representatives of $\Sym_\mu$ in $\Sym_n$. Then, as an $R$-module, 
$$
\Ind_{\H(\Sym_\mu)}^\H(V)\cong\bigoplus_{d\in\D_\mu}V\otimes T_d
$$
by \cite[Theorem~2.7]{DJ:reps}.

Let $x_\mu = \sum_{w\in\Sym_\mu}T_w$. Then $T_w x_\mu = x_\mu T_w =
q^{\len(w)}x_\mu$ for all $w\in\Sym_\mu$.  The {\sf trivial
representation} of $\H(\Sym_\mu)$ is the free $R$--module 
$\mathbf 1_\mu = Rx_\mu$. 

Let $y_\mu = \sum_{w\in\Sym_\mu} (-q)^{-\len(w)} T_w$. Then 
$T_w y_\mu = y_\mu T_w = (-1)^{\len(w)} y_\mu$ for all $w\in\Sym_\mu$. 
The {\sf sign representation} of $\H(\Sym_\mu)$ is the free $R$--module 
$\mathcal E_\mu = Ry_\mu$.

For any composition $\mu$ we define the permutation module
$M(\mu)=\Ind_{\H(\Sym_\mu)}^\H(\mathbf 1_\mu)\cong x_\mu\H$. Then
$M(\mu)$ is free as an $R$--module of rank $[\Sym_n:\Sym_\mu]$ with
basis $\set{x_\mu T_d|d\in\D_\mu}$. The $\H$--action on $M(\mu)$ is
determined by 
$$x_\mu T_dT_i=\begin{cases}
          qx_\mu T_d,&\If \len(dr_i)>\len(d)\And dr_i\not\in\D_\mu,\\
 \phantom q x_\mu T_{dr_i},&\If\len(dr_i)>\len(d)\And dr_i\in\D_\mu,\\
          qx_\mu T_{dr_i}+(q-1)x_\mu T_d,&\Otherwise.
\end{cases}$$
Note that if $\len(dr_i)<\len(d)$ then $dr_i\in\D_\mu$. 

Let $*\map\H\H$ be the $R$--linear map on $\H$ determined by
$T_w^* = T_{w^{-1}}$, for all $w\in\Sym_n$. This defines
an $R$-algebra anti--automorphism on $\H$ of order $2$.

The module $M(\mu)$ carries a symmetric bilinear form 
$\<\ ,\ \>_\mu$ given by 
$$
\<x_\mu T_a,x_\mu T_b\>_\mu \; =\;
   \begin{cases}
                  q^{\len(a)}, & \If a=b,\\
                  0,           & \Otherwise,
   \end{cases}
$$
for $a,b\in\D_\mu$. It follows from the formulae above that the form
$\<\ ,\ \>_\mu$ is associative in the sense that
$$
\<xh,y\>_\mu\; =\;\<x,y h^*\>_\mu\;
$$
for all $x,y\in M(\mu)$ and all $h\in\H$. 

We will need two dualities on the category of right $\H$--modules.
Both of them come from involutions on $\H$. The first duality comes
from the involution $*$ defined above. The second is induced from the
automorphism $\#\map\H\H$ which is the $R$--linear map on $\H$
determined by $T_w^\# = (-q)^{l(w)}T_{w^{-1}}^{-1}$, for all
$w\in\Sym_n$. It is straightforward to check that $\#$ preserves the
relations in $\H$ and, hence, that it is an $R$-algebra automorphism of order $2$.
Note that the involutions $\#$ and $*$ commute.

If $V$ is an $\H$-module let $V^*$ be its $R$-linear dual. Then $V^*$
becomes an $\H$--module by letting $(\phi\cdot\xi)(v) :=
\phi(v\xi^*)$, where $\phi\in V^*$, $v\in V$ and $\xi\in\H$. With the
according operation on morphisms, this
defines a contravariant self--equivalence on the category of
$\H$-modules.

If $V$ is an $\H$-module let $V^\#$ the $\H$-module with underlying
$R$-module $V$ and operation $v\cdot_\# \xi := v\cdot \xi^\#$, where $v\in V$ and
$\xi\in\H$. With the identical operation on morphisms, this defines a covariant 
self--equivalence on the category of $\H$-modules.

\section{Specht modules} 

We recall some well-known facts due to Dipper and
James~\cite{DJ:reps}.

Let $\lambda = (\lambda_1,\lambda_2,\dots)$ be a composition of $n$.
The {\sf diagram} of $\lambda$ is the set $[\lambda] = \set{(i,j)\in\N^2|1\le
j\le\lambda_i}$. We identify the diagram of $\lambda$ with an array of
boxes in the plane. For example, if $\lambda = (4,3,2)$ then
$$
[\lambda] \; =\; \yng(4,3,2) \; .
$$
The {\sf conjugate} of $\lambda$ is the partition
$\lambda' = (\lambda_1',\lambda_2',\dots)$, where
$\lambda_j' = \#\set{i\ge 1|\lambda_i\ge j}$ for all $j$; that is,
$\lambda'$ is the partition of $n$ whose diagram is obtained by
interchanging the rows and columns of the diagram of $\lambda$.

Formally, a {\sf $\lambda$--tableau} is a bijection
$\t\map{[\lambda]}\{1,2,\dots,n\}$; however, we will think of a
$\lambda$--tableau as a labelling of the diagram of $\lambda$ by the
numbers $1,2,\dots,n$. Accordingly, we will speak of the rows and
columns of a tableau. For example, 
$$\young(123,45)\;,\qquad\young(135,24)\;,\qquad\young(145,23)\qquad\And\qquad
  \young(234,15)$$ 
are all $(3,2)$--tableaux. 

A tableau is {\sf row standard} if in each row its entries increase
from left to right. A tableau is {\sf standard} if it is row standard
and in each column its entries increase from top to bottom. Let
$\Std(\lambda)$ be the set of standard $\lambda$-tableaux.

All of the tableaux above are row standard; however, only the first
two tableaux are standard.

The {\sf initial $\lambda$--tableau} $\t^\lambda$ is the standard
$\lambda$--tableau which has the numbers $1,2,\dots,n$
entered in order from left to right, and then top to bottom, along its
rows. The {\sf terminal $\lambda$--tableau} $\t_\lambda$ is the standard
$\lambda$--tableau which has the numbers $1,2,\dots,n$
entered in order from top to bottom, and then left to right, along its
columns. Of the $(3,2)$--tableaux above, the first is $\t^{(3,2)}$ and
the second is $\t_{(3,2)}$.

The symmetric group $\Sym_n$ acts from the right on the set of
$\lambda$--tableaux by permuting their entries. If $\t$ is a
$\lambda$--tableau let $d(\t)$ be the unique permutation such that
$\t=\t^\lambda d(\t)$. In particular, we set
$w_\lambda=d(\t_\lambda)$.

We remark that $\D_\mu = \set{d(\t)|\t\text{\ is a row standard $\mu$--tableau}}.$

Suppose that $\lambda$ is a partition of $n$ and let 
$z_\lambda= x_\lambda T_{w_\lambda} y_{\lambda'}$. The {\sf Specht
module} is the submodule $S(\lambda) = z_\lambda\H$ of $M(\lambda)$.

Let $S(\lambda)^\perp = \set{x\in M(\lambda)|\< x,y\>=0\ForAll y\in S(\lambda)}$.
As $\<\ ,\ \>_\lambda$ is associative,
$S(\lambda)^\perp$ is an $\H$-submodule of $M(\lambda)$. More
precisely, $S(\lambda)^\perp$ is the kernel of the $\H$-linear map
$$
M(\lambda)\;\lraa{\delta_\lambda}\;S(\lambda)^\ast \; ;\;\;\; x_\lambda h \;\lramaps\;\< x_\lambda h, - \>_\lambda \; ,
$$
where $h\in\H$. 

By restricting the bilinear form $\<\ ,\ \>_\lambda$ on $M(\lambda)$
we obtain a bilinear form on $S(\lambda)$.  If $R$ is a field then
$D(\lambda) = S(\lambda)/S(\lambda)\cap S(\lambda)^\perp$ is either
zero or absolutely irreducible. Moreover, all of the irreducible
$\H$--modules arise uniquely in this way \cite[Theorem~5.2]{DJ:reps}.

Before we can give a basis of $S(\lambda)$ we need some more notation.
If $\t$ is a $\lambda$--tableau let $\t'$ be the $\lambda'$--tableau
obtained by interchanging the rows and columns of~$\t$. For example,
$(\t^\lambda)'=\t_{\lambda'}$ and $(\t_\lambda)'=\t^{\lambda'}$.
Finally, if $\t$ is a standard $\lambda$--tableau let 
$v_\t = z_\lambda T_{d(\t')}$.

\begin{Point}\it{(Dipper--James~\cite[Theorem~5.6]{DJ:reps})}
The Specht module $S(\lambda)$ is free as an $R$--module with basis 
$\set{v_\t|\t\in\Std(\lambda)}$.
\label{DJ_basis}\end{Point}

We call $\set{v_\t|\t\in\Std(\lambda)}$ the {\sf Dipper--James basis} of $S(\lambda)$.
Let $n_\lambda = \#\Std(\lambda)$ be the number of standard
$\lambda$--tableaux. Then, as an $R$-module, $S(\lambda)$ is free of
rank~$n_\lambda$.

Fix an ordering of $\Std(\lambda)$ and let
$$G(\lambda)=\Big(\<v_\s,v_\t\>_\lambda\Big)_{\s,\t\in\Std(\lambda)}$$ 
be the Gram matrix of the bilinear form $\<\ ,\ \>_\lambda$, with
respect to the Dipper--James basis. The matrix $G(\lambda)$ depends on
the choice of ordering on $\Std(\lambda)$; however, all of the quantities
that we are interested in will be independent of this choice. We remark
that $\det G(\lambda)$ has been explicitly computed by Dipper and
James~\cite[Theorem~4.11]{DJ:blocks}.

\section{Diagonalizability and elementary divisors}

Given an integer $m\geq 1$, an $m\times m$ matrix $A$ with coefficients in $R$ is
{\sf diagonalizable} if there exist matrices $S$ and $T$ in $\GL_m(R)$ such that
$SAT$ is a diagonal matrix. The matrix $A$ is {\sf divisibly diagonalizable}
if $SAT = \diag(d_1,\dots,d_m)$ 
is a diagonal matrix such that $d_i$ divides $d_{i+1}$ in $R$, 
for $1\leq i < m$. If $A$ is divisibly diagonalizable and
$SAT = \diag(d_1,\dots,d_m)$ satisfies this condition, then we call $d_1,\dots,d_m$ 
the elementary divisors of $A$.

Given $A\in R^{m\times m}$, we let $I_k(A)$ be the ideal of the
$k\times k$ minors of $A$, for $1\leq k\leq m$. Note that for $B\in
R^{m\times m}$, we have $I_k(AB) \ss I_k(A)$ and $I_k(BA) \ss I_k(A)$.
Hence for $S,T\in\GL_m(R)$, we have $I_k(A) = I_k(SAT)$. Therefore, if
$A$ is divisibly diagonalizable with resulting diagonal elements
$d_1,\dots,d_m$, then $I_k(A) = I_k(\diag(d_1,\dots,d_m))$ is the
principal ideal generated by $d_1 d_2\cdots d_k$. This shows that the
resulting diagonal entries are independent, up to multiplication by
units, of the choice of the diagonalizing matrices.  In other words,
the elementary divisors of a divisibly diagonalizable matrix are
well--defined modulo units.

Whether or not $A$ is diagonalizable, the ideals $I_k(A)\ss R$ are invariant
under the equivalence relation $A\sim SAT$. It would be interesting to consider
the equivalence classes within 
$\set{A\in R^{m\times m} | I_k(A) = J_k\For 1\le k\le m}$ for 
a fixed tuple $(J_i)$ of ideals of $R$.

If $R$ is a principal ideal domain then every matrix $A\in R^{m\times m}$ is divisibly 
diagonalizable by the elementary divisor theorem.
The resulting diagonal matrix is known as the Smith normal form.  

Now the Laurent polynomial ring $\Z[q,q^{-1}]$ is not a principal
ideal domain and, in fact, there are strict inclusions of the set of
divisibly diagonalizable matrices in the set of diagonalizable
matrices, and of the set of diagonalizable matrices in all matrices
with coefficients in $\Z[q,q^{-1}]$. For example, the matrix $A =
\big(\begin{smallmatrix} q-1 & 0 \\ 0 & q+1 \end{smallmatrix}\big)$ is
diagonalizable, but not divisibly diagonalizable because $I_1(A)$, the
ideal of $R$ generated by the entries of $A$, is not principal. 

Proving that a matrix is not diagonalizable is slightly harder. For example, 
we claim that the matrix
$B = \big(\begin{smallmatrix} q+1 & 2 \\ 0 & q+1 \end{smallmatrix}\big)$
    is not diagonalizable over $\Z[q,q^{-1}]$. To see this, notice that over $\Q[q,q^{-1}]$
the matrix $B$ has elementary divisors $1$ and $(q+1)^2$. Therefore,
if $B$ is diagonalizable over $\Z[q,q^{-1}]$ then one of these
diagonal entries must be a unit in $\Q[q,q^{-1}]$; that is, of the
form $aq^b$ with $a,b\in\Z$.  Reducing modulo $2$ this shows that one
of the elementary divisors of $B$ over $\F_2[q,q^{-1}]$ is zero or a
unit. However, this is a contradiction because the elementary
divisors of $B$ over $\F_2[q,q^{-1}]$ are $q+1$ and $q+1$.

In proving that certain Gram matrices $G(\lambda)$ are  divisibly diagonalizable over 
$\Z[q,q^{-1}]$, we shall make use of the following simple lemma. 

\begin{Lemma}\label{diag}
Let $A$ be an $m\times m$ matrix with coefficients in $R$, 
and suppose that there exist invertible matrices $S,T\in\GL_m(R)$ such that 
$$
SAT = \left(\begin{matrix} d_1    & b_{12} & \dots  & b_{1m} \\
                           0      & d_2    & \dots  & b_{2m} \\
                           \vdots & \ddots & \ddots & \vdots\\
                           0      & \cdots & 0      & d_m
\end{matrix}\right),
$$
where $d_1\, |\, d_2\, |\, \dots\, |\, d_m$ and $d_i$ divides $b_{ij}$ for all $j$. Then
$A$ is divisibly diagonalizable and $d_1,d_2,\dots,d_m$ are the
elementary divisors of $A$.
\end{Lemma}

\begin{proof}
The matrix $SAT$ can be written as the product of $\diag(d_1,\dots,d_m)$ with a 
matrix in $\GL_m(R)$.
\end{proof}

As we saw with the non--diagonalizable matrix 
$\big(\begin{smallmatrix} q+1 & 2 \\ 0 & q+1 \end{smallmatrix}\big)$
above, the requirement that $d_i$ divides $b_{ij}$ for all $j$ is not
superfluous. 

\section{Elementary divisors for conjugate partitions}
\label{SecConj}

Let $R = \Z[q,q^{-1}]$. Let $\lambda$ be a partition of $n$. In this
section we relate the Gram matrices $G(\lambda)$ and $G(\lambda')$. We
start with some mild generalizations of some results about Specht
modules which were proved by Dipper and James~\cite{DJ:reps} over a
field.

Recall that if $Y$ is a submodule of an $R$--free module $X$ then $Y$
is \textsf{pure} if the quotient module $X/Y$ is $R$-free.

\begin{Lemma}
\label{LemTr1}
Suppose that $\lambda$ is a partition. Then the Specht module
$S(\lambda)$ is a pure submodule of $M(\lambda)$. 
\end{Lemma}

\begin{proof}
Using the Dipper--James basis of $S(\lambda)$, and the basis
$\set{x_\lambda T_d|d\in\D_\lambda}$ of $M(\lambda)$, suitably ordered,
the matrix representing the embedding $S(\lambda)\lra M(\lambda)$
$\Z[q,q^{-1}]$-linearly becomes triangular with $1$s on the diagonal
\cite[Theorem~5.8]{DJ:reps}. 
\end{proof}

\begin{Corollary}
The map $M(\lambda)\lraa{\delta_\lambda}S(\lambda)^*$ is surjective.
\end{Corollary}

\begin{proof}
The map $\delta_\lambda$ is the composition of the map
$$M(\lambda)\; \lraiso \; M(\lambda)^*\; ;\;\;\; \xi\;\lramaps \;\<\xi, -\>_\lambda$$ 
with the dual of the inclusion map $S(\lambda)\lra M(\lambda)$. This
is surjective by (\ref{LemTr1}).
\end{proof}

Let $\alpha(\lambda) = \sum_{i\geq 1} (i-1)\lambda_i 
                     = \sum_{i\geq 1} \binom{\lambda'_i}{2}$. 
Note that $\alpha(\lambda)=\len(w_{0,\lambda'})$, where
$w_{0,\lambda'}$ is the unique length of longest element in
$\Sym_{\lambda'}$. The next lemma is well known; see, for example,
\cite[Prop.~2.2]{Takeuchi}. We include a proof for completeness.

Recall that automorphism $\#$, and the corresponding operation on the
module category of $\H$, were defined at the end of section~2.

\begin{Lemma}\label{hash}
We have $x_\lambda^\# = q^{\alpha(\lambda')} y_\lambda$ and 
$y_\lambda^\# = q^{-\alpha(\lambda')} x_\lambda$.
\end{Lemma}

\begin{proof}
As $\#$ is an involution the two equalities are equivalent, so we
prove only the first.  For any integer $i$, with $1\le i < n$, we have
$x_\lambda^\# T_i=(x_\lambda T_i^\#)^\# = -x_\lambda^\#$. Write
$x_\lambda^\# = \sum_{w\in \Sym_\lambda} a_w T_w$, for some
$a_w\in\Z[q,q^{-1}]$.  Comparing coefficients on both sides of the
equation $x_\lambda^\# T_i=-x_\lambda^\#$ shows that $a_{wr_i} = (-q)
a_w$ for each $w$ that has a reduced expression ending in $r_i$;
compare \cite[Cor.~1.7]{M:IHA}. Hence, $x_\lambda^\#$ is a scalar multiple of
$y_\lambda$. Then
$T_{w_{0,\lambda}}^\#=(-1)^{\len(w_{0,\lambda})}T_{w_{0,\lambda}}$
plus a linear combination of $T_v$ where $v\in\Sym_\lambda$ and
$\len(v)<\len(w_{0,\lambda})$. Therefore, comparing the coefficient of
$T_{w_{0,\lambda}}$ in $x_\lambda^\#$ and $y_\lambda$ gives the
result.
\end{proof}

Recall that $S(\lambda)=z_\lambda\H$, where 
$z_\lambda=x_\lambda T_{w_\lambda} y_{\lambda'}$. The importance of
$z_\lambda$, and the irreducibility of $S(\lambda)$ in the semisimple
case, follow from the following simple fact.

\begin{Point}\it{(Dipper--James~\cite[Lemma~4.1]{DJ:reps})}\label{vanishing}
Suppose that $w\in\Sym_n$. Then 
$$x_\lambda T_w y_{\lambda'}=\begin{cases}
       \pm q^a z_\lambda,&\If w\in\Sym_\lambda w_\lambda\Sym_{\lambda'},\\
                       0,&\Otherwise,
\end{cases}$$
for some integer $a$.
\end{Point}

The proof of this result amounts to the observation that
$\Sym_\lambda\cap w\Sym_{\lambda'}w^{-1}=\{1\}$ if and only if
$w\in\Sym_\lambda w_\lambda\Sym_{\lambda'}$. 

\begin{Lemma}[The Submodule Theorem]
\label{LemJST}
If $U$ is a pure submodule of $M(\lambda)$, then $S(\lambda)\ss U$ or $U\ss S(\lambda)^\perp$.
\end{Lemma}

\begin{proof}
For all $u\in U$, we have $u y_{\lambda'} = \alpha_u z_\lambda$ for some $\alpha_u\in\Z[q,q^{-1}]$ by (\ref{vanishing}). 

{\it Case $1$}: $\alpha_u = 0$ for all $u\in U$. Therefore, if $u\in U$ and
$h\in\H$ then we have $\<u,z_\lambda h\>_\lambda = \<u h^*
y_{\lambda'},x_\lambda T_{w_\lambda}\>_\lambda$, since $y_{\lambda'}^*
= y_{\lambda'}$. But $u h^*\in U$, so $u h^* y_{\lambda'}
= 0$ and $u\in S(\lambda)^\perp$. Hence, $U\ss S(\lambda)^\perp$.

{\it Case $2$}: $\alpha_u \neq 0$ for some $u\in U$. Now $U\ni u
y_{\lambda'} = \alpha_u z_\lambda$ implies $z_\lambda\in U$ since
$U\ss M^\lambda$ is a pure submodule. Therefore, $S(\lambda)\ss U$.
\end{proof}

Note that the right ideal $y_{\lambda'} T_{w_\lambda}^{-1} x_\lambda\H$ is isomorphic to $S(\lambda)^\#$ via $\xi\lramaps\xi^\#$. Composing left multiplication
by $y_{\lambda'} T_{w_\lambda}^{-1}$ with this isomorphism, and using
Lemma~\ref{hash}, we obtain a surjective $\H$-linear map
$$
M(\lambda)\; =\; x_\lambda \H \; \lraa{\theta_\lambda} \; S(\lambda')^\# \; ; \;\; x_\lambda h \; \lramaps \; z_{\lambda'}\,\cdot_\#\, h\; , 
$$
where $h\in\H$.

\begin{Lemma}
We have $\Kern\theta_\lambda = S(\lambda)^\perp$. This induces an isomorphism
$$
S(\lambda')^\#\;\lraisoa{\psi_\lambda}\;S(\lambda)^* \; ;\;\;\; z_{\lambda'}\,\cdot_\#\, h \;\lramaps\;\< x_\lambda h , -\>_\lambda \; ,
$$
where $h\in\H$.
\end{Lemma}

\begin{proof}
Both $\Kern\theta_\lambda$ and $S(\lambda)^\perp$ are pure submodules
of $M(\lambda)$.  Over $\Z[q,q^{-1}]$, both $S(\lambda')^\#$ and
$S(\lambda)^*$ are free of rank $n_\lambda$, so it suffices to prove
that $\Kern\theta_\lambda\ss S(\lambda)^\perp$. By (\ref{LemJST}) this
is equivalent to showing that $S(\lambda)\not\ss \Kern\theta_\lambda$.
So it is enough to show that $z_\lambda \theta_\lambda \neq 0$. 
The bilinear form $\<\ ,\ \>_{\lambda'}$ is associative, so
\begin{align*}
\<z_\lambda \theta_\lambda, x_{\lambda'}\>_{\lambda'}
& =  q^{-\alpha(\lambda)}\< z_{\lambda'} T_{w_{\lambda'}}^{-1} x_{\lambda'} , x_{\lambda'}\>_{\lambda'} \\
& =  q^{-\alpha(\lambda')}(\Sum_{w\in\Sym_{\lambda'}} q^{\len(w)})\<z_{\lambda'}, x_{\lambda'} T_{w_\lambda}^{-1}\>_{\lambda'}\; . \\
\end{align*}
Now, 
$z_{\lambda'} = \sum_{v\in\Sym_\lambda}(-q)^{-\len(v)} x_{\lambda'} T_{w_{\lambda'} v}$, where each $w_{\lambda'} v$ is a distinguished coset representative for 
$\Sym_{\lambda'}$. In contrast, $T_{w_\lambda}^{-1}$ is equal to
$T_{w_{\lambda'}}$ plus a linear combination of terms $T_u$, where
$u\in\Sym_n$ with $\len(u) < \len(w_{\lambda'})$. Thus
$\<z_{\lambda'}, x_{\lambda'} T_{w_\lambda}^{-1}\>_{\lambda'} 
       =\<x_{\lambda'} T_{w_{\lambda'}}, x_{\lambda'} T_{w_{\lambda'}}\>_\lambda 
       = q^{\len(w_{\lambda'})}$. 
Hence, $\<z_\lambda \theta_\lambda, x_{\lambda'}\>_{\lambda'} \neq 0$.

A comparison of the short exact sequences
$$
0\;\lra\;\Kern\theta_\lambda\;\lra\; M(\lambda) \;\lraa{\theta_\lambda}\; S(\lambda')^\# \;\lra\; 0 
$$
and 
$$
0\;\lra\; S(\lambda)^\perp\;\lra\; M(\lambda) \;\lraa{\delta_\lambda}\; S(\lambda)^*\;\lra\; 0
$$
yields the isomorphism $\psi_\lambda$.
\end{proof}

For each node $(i,j)\in [\lambda]$, we let 
$h_{i,j} = (\lambda_i - j) + (\lambda'_j - i)+1$ be the corresponding
hook length and set 
$h_{\lambda}(q) = \prod_{(i,j)\in [\lambda]} [h_{i,j}]_q$. The next
lemma follows from results of {Murphy} \cite{murphy:hecke}.

\begin{Lemma}\label{scalar}
We have $z_\lambda T_{w_\lambda}^{-1} z_\lambda = q^{n - \alpha(\lambda)} h_\lambda(q) z_\lambda$. 
\end{Lemma}

\begin{proof} For the purpose of this proof, we may assume $R = \Q(q)$. By 
\cite[p.~510--511]{murphy:hecke}, there exists an element 
$\Psi^*_{\t_\lambda}=T_{w_{\lambda'}}+\sum_{\len(v)<\len(w_{\lambda'})}r_v T_v\in\H$, 
for some $r_v\in R$, such that
\[
z_\lambda\Psi^*_{\t_\lambda}\; =\; q^{n - \alpha(\lambda) + \len(w_{\lambda'})} h_\lambda(q) E_\lambda\; ,
\]
where $E_\lambda$ is a primitive idempotent such that
$E_\lambda \H = z_\lambda\H = S(\lambda)$.  In particular 
$E_\lambda z_\lambda=z_\lambda$. (Note that $z_\lambda=z_{\lambda t}$ 
in Murphy's notation; see \cite[p.~496, p.~498]{murphy:hecke}.)

Note that $T_{w_\lambda}^{-1} = q^{-\len(w_\lambda)}T_{w_{\lambda'}} +
\sum_{\len(v) < \len(w_{\lambda'})} r_v T_v$, for some $r_v\in R$.
Now, if $\len(v)<\len(w_{\lambda'})$ then 
$v \not\in \Sym_{\lambda'} w_{\lambda'} \Sym_\lambda$, so
$y_{\lambda'} T_v x_\lambda =(x_{\lambda} T_v^* y_{\lambda'})^*= 0$ by
(\ref{vanishing}). Consequently,  $z_\lambda T_v z_\lambda = 0$.
Therefore,
\begin{xalignat*}{2}
z_\lambda T_{w_\lambda}^{-1}z_\lambda
   &= q^{-\len(w_\lambda)} z_\lambda T_{w_{\lambda'}} z_\lambda \\
   &= q^{-\len(w_\lambda)} z_\lambda \Psi^*_{\t_\lambda} z_\lambda \\
   &= q^{n - \alpha(\lambda)} h_\lambda(q) E_\lambda z_\lambda \\
   &= q^{n - \alpha(\lambda)} h_\lambda(q)  z_\lambda\; .
\end{xalignat*}
  
\end{proof}

Consider the $\H$-linear map
\[
S(\lambda)\;\lraa{\gamma_\lambda}\;S(\lambda)^\ast \; ;\;\;\; \xi\;\lramaps\;\< \xi , -\>_\lambda\; .
\]

\begin{Lemma}
\label{PropComp}
The composition
\[
S(\lambda)\;\lraa{\gamma_\lambda}\; S(\lambda)^\ast \;\lraisoa{\psi_\lambda^{-1}}\; S(\lambda')^\# \;\lraa{\gamma_{\lambda'}^\#}\; S(\lambda')^{\ast,\#} 
\;\;\lraisoa{(\psi_{\lambda'}^\#)^{-1}}\;\; S(\lambda)
\]
is equal to scalar multiplication by $(-q)^{\len(w_{\lambda'})} q^{n - \alpha(\lambda) - \alpha(\lambda')} h_\lambda(q)$.
\end{Lemma}

\begin{proof}
The element $z_\lambda$ is mapped via $\gamma_\lambda$ to $\< z_\lambda , -\>_\lambda$, which is mapped via 
$\psi_\lambda^{-1}$ to $z_{\lambda'}\,\cdot_\#\, T_{w_\lambda} y_{\lambda'} = z_{\lambda'} T_{w_\lambda}^\# y_{\lambda'}^\#$, which in turn
goes to $\< z_{\lambda'} T_{w_\lambda}^\# y_{\lambda'}^\# , -\>_\lambda$ via $\gamma_{\lambda'}^\#$, and finally to 
\[
\begin{array}{rcl}
z_\lambda\,\cdot_\#\, T_{w_{\lambda'}} y_\lambda T_{w_\lambda}^\# y_{\lambda'}^\#  
& = & (-q)^{\len(w_{\lambda'})} q^{-\alpha(\lambda')} z_\lambda T_{w_\lambda}^{-1} z_\lambda \vspace{1mm}\\
& = & (-q)^{\len(w_{\lambda'})} q^{n - \alpha(\lambda) - \alpha(\lambda')} h_\lambda(q) z_\lambda \\
\end{array}
\]
via $(\psi_{\lambda'}^\#)^{-1}$, by Lemma~\ref{scalar}.
\end{proof}

Let $I_m$ be the $m\times m$ identity matrix. Recall that
$n_\lambda=\#\Std(\lambda)$ is the dimension of the Specht module
$S(\lambda)$.

\begin{Proposition} 
\label{CorDiagbar}
Suppose that $\lambda$ is a partition of $n$.
\\\vspace*{-4mm}
\begin{itemize}
\item[(1)] There exist invertible matrices $A,B\in\GL_{n_\lambda}(\Z[q,q^{-1}])$ 
such that 
\[
G(\lambda)\cdot A \cdot G(\lambda')\cdot  B \; =\; h_\lambda(q)\cdot I_{n_{\lambda}}\; .
\] 
\item[(2)] $G(\lambda)$ is diagonalizable to the diagonal matrix $D$ if and only if $G(\lambda')$ is diagonalizable to the diagonal matrix $h_\lambda(q) D^{-1}$.
\item[(3)] $G(\lambda)$ is divisibly diagonalizable if and only if $G(\lambda')$ is divisibly diagonalizable. In this case, the product of the $i^\th$ 
elementary divisor of $G(\lambda)$ and the $(n_\lambda + 1 - i)^\th$ elementary divisor of $G(\lambda')$ is equal to $h_\lambda(q)$.
\end{itemize}
\end{Proposition}

Recall that elementary divisors are only well defined up to a unit in
$\Z[q,q^{-1}]$; the same is true of their product in (3).

\begin{proof}
(1) The $R$-linear map $\gamma_\lambda$ is represented by the matrix $G(\lambda)$ with respect to the Dipper--James basis and its dual basis. Thus
the assertion follows by (\ref{PropComp}).

(2) If $G(\lambda) = SDT$ with $S,T\in\GL_{n_\lambda}(\Z[q,q^{-1}])$ and $D\in\Z[q,q^{-1}]^{n_\lambda\times n_\lambda}$ is a diagonal matrix, then 
$G(\lambda') = A^{-1}T^{-1}(h_\lambda(q) D^{-1})S^{-1}B^{-1}$. Since $G(\lambda')$ has coefficients in $\Z[q,q^{-1}]$, so does $h_\lambda(q) D^{-1}$.

(3) Repeat the argument of (2). 
\end{proof}

We remark that all of the results in this section hold more generally
when the Hecke algebra $\H$ is defined over an integral domain $R$
such that $\H\otimes_R Q$ is semisimple, where $Q$ is the field of
fractions of $R$. (We need semisimplicity over~$Q$ only when we apply
Murphy's results in the proof of Lemma~\ref{scalar}.) In particular,
Proposition~\ref{CorDiagbar} holds when $R=F[q,q^{-1}]$ and $F$ is any
field. Notice that $G(\lambda)$ is always diagonalizable in this case
because $F[q,q^{-1}]$ is a principle ideal domain.

\section{The elementary divisors for hook partitions}

Throughout this section we fix an integer $k$, with $0\le k<n$, and
consider the Specht module $S(\lambda)$, where $\lambda=(n-k,1^k)$. We
will show that $G(\lambda)$ is divisibly diagonalizable by explicitly
constructing two bases of $S(\lambda)$ which transform $G(\lambda)$
into the form of the matrix in Lemma~\ref{diag}. In particular, this
will allow us to determine the elementary divisors of $S(\lambda)$.

The Specht module $S(\lambda)$ is defined as a submodule of the
permutation module $M(\lambda)$; however, to compute the elementary
divisors we will work inside a different permutation module.

By definition,
$S(\lambda)=x_\lambda T_{w_\lambda}y_{\lambda'}\H
           =x_{(n-k,1^k)}T_{w_{(n-k,1^k)}}y_{(k+1,1^{n-k-1})}$.
We first need to understand the permutation $w_\lambda=w_{(n-k,1^k)}$
a little better. This requires some new notation. For integers
non--negative $i$ and $j$ define
$$r_{i,j}=\begin{cases}1,                   &\If i=0\Or j=0,\\
                       r_ir_{i+1}\dots r_j, &\If 0<i\le j,\\
                       r_ir_{i-1}\dots r_j, &\If i>j>0,
\end{cases}$$
and set $T_{i,j}=T_{r_{i,j}}$.  Next, let $a$ and $b$ be non--negative
integers such $a+b\le n$. If either $a=0$ or $b=0$ then set
$w_{a,b}=1$. If both $a$ and $b$ are non--zero then define
$w_{a,b}=(r_{a+b-1,1})^b$; then one can check that, in two--line
notation,
$$w_{a,b}=\Big(\begin{array}{cccc|cccc}
                 1&2&\dots&a&a+1&a+2&\dots&a+b\\
              b+1&b+2&\dots&a+b&1&2&\dots&b
\end{array}\Big).$$
It is not hard to see that $w_{b,a}=w_{a,b}^{-1}$
and that $w_{a,b}=r_{a,a+b-1}w_{a-1,b}$ and 
$\len(w_{a,b})=\len(r_{a,a+b-1})+\len(w_{a-1,b})$; see \cite{DJ:B}. 
Consequently,
$$w_{a,b}=r_{a,a+b-1}r_{a-1,a+b-2}\dots r_{1,b}
         =r_{1,a}r_{2,a+2}\dots r_{b,a+b-1}$$
with the lengths adding in both cases. Hence, $\len(w_{a,b})=ab$.

The permutation $w_{(n-k,1^k)}$ is essentially one of these
permutations because
$$w_{(n-k,1^k)}=\Big(\begin{array}[c]{*7c}
           1&2&\cdots&n-k&n-k+1&\cdots&n\\
           1&k+2&\cdots&n&2&\cdots&k+1
\end{array}\Big).$$
Hence, $w_{(n-k,1^k)}=r_{n-k,n-1}r_{n-k-1,n-2}\dots r_{2,k+1}$, with
the lengths adding. So
$$T_{w_{(n-k,1^k)}}=T_{n-k,n-1}\dots T_{2,k+1},$$ 
Notice also that $T_{w_{n-k.k}}=T_{w_{(n-k,1^k)}}T_{1,k}$.

If $w\in\Sym_{(k,n-k)}\cong\Sym_k\times\Sym_{n-k}$
then we write $w=(u,v)$, where $u\in\Sym_k$ and $v\in\Sym_{(1^k,n-k)}$
are the unique permutations such that $w=uv=vu$. Set 
$$x_{(k|n{-}k)}=y_{(k,1^{n-k})}x_{(1^k,n-k)}
       =\sum_{(u,v)\in\Sym_{(k,n-k)}}(-q)^{-\len(u)}T_{uv}.$$
Then it is easy to see that $Rx_{(k|n-k)}$ is an $\H(\Sym_\mu)$--module on
which the subalgebras $\H(\Sym_k)$ and $\H(\Sym_{(1^k,n-k)})$ act via
their sign and trivial representations, respectively. Let
$$M(k|n{-}k)=\Ind_{\H(\Sym_{(k,n-k)})}^{\H(\Sym_n)}\big(Rx_{(k|n{-}k)}\big)
              \cong x_{(k|n{-}k)}\H.$$ 
As in section~2, the induced module $M(k|n{-}k)$ is free as an $R$--module with
basis $\set{x_{(k{\mid}n-k)}T_d|d\in\D_{(k,n-k)}}$. Furthermore,
$M(k|n{-}k)$ possesses a natural non--degenerate associative bilinear form 
$\<\ ,\ \>_{(k|n{-}k)}$ which is determined by 
$$\<x_{(k|n{-}k)} T_u,x_{(k|n{-}k)} T_v\>_{(k|n{-}k)}
  =\begin{cases}
                  q^{\len(u)}, & \If u=v,\\
                  0,&\Otherwise,
   \end{cases}
$$
for $u,v\in\D_{(k,n-k)}$. Donkin~\cite{D:super} calls $M(k|n{-}k)$ a
{\sf trivial source module}.

Let $y_{k+1}'=1+\sum_{j=1}^k(-q)^{j-k-1}T_{k,j}
        =1-q^{-1}T_k+q^{-2}T_{k,k-1} +\dots+(-q)^{-k}T_{k,1}$. This is
a sum over the right coset representatives of $\Sym_k$ in
$\Sym_{k+1}$. Consequently, it follows that
$y_{(k+1,1^{n-k-1})}=y_{(k,1^{n-k})}y_{k+1}'$.
The reason for introducing the module $M(k|n{-}k)$ is the following
result.

Given a non--negative integer $k>1$ let $[k]_q =1+q+\dots+q^{k-1}$ and
$[k]_q^! = [1]_q [2]_q \cdots [k]_q$. Notice that if $q=1$ then
$[k]_1=k$ and $[k]_1^!=k!$.

\begin{Proposition}
Let $\lambda = (n-k,1^k)$. The map 
$$\pi_k\;\map{\; S(\lambda)\;}\; M(k|n{-}k)\; ;\;\; 
    z_\lambda h\;\mapsto\; x_{(k|n{-}k)}y_{k+1}'h$$
is an injective $\H$--module homomorphism. Moreover, 
$$\<x,y\>_\lambda 
   =q^{\frac k2(2n-3k-1)}[k]_q^!\<\pi(x),\pi(y)\>_{(k|n{-}k)},$$
for all $x,y\in S(\lambda)$.
\label{changing}
\end{Proposition}

\begin{proof}
By definition, 
$S(\lambda)=x_{(n-k,1^k)}T_{w_\lambda} y_{(k+1,1^{n-k-1})}\H$.
As remarked above, $w_{n-k,k}=w_{(n-k,1^k)}r_{1,k}$ with
the lengths adding. Therefore, since $r_{1,k}\in\Sym_{\lambda'}$, 
\begin{align*}
x_{(n-k,1^k)}T_{w_{(n-k,1^k)}} y_{(k+1,1^{n-k-1})}
          &=(-1)^kx_{(n-k,1^k)}T_{w_{(n-k,1^k)}}T_{1,k}y_{(k+1,1^{n-k-1})}\\
          &=(-1)^kx_{(n-k,1^k)}T_{w_{n-k,k}}y_{(k+1,1^{n-k-1})}\\
          &=(-1)^kT_{w_{n-k,k}}x_{(1^k,n-k)}y_{(k+1,1^{n-k-1})}\\
          &=(-1)^kT_{w_{n-k,k}}x_{(1^k,n-k)}y_{(k,1^{n-k})}y_{k+1}'\\
          &=(-1)^kT_{w_{n-k,k}}x_{(k|n{-}k)}y_{k+1}'.
\end{align*}
Therefore, $\pi(x)=(-1)^kT_{w_{n-k,k}}^{-1}x$, for all $x\in S(\lambda)$. As
$T_{w_{n-k,k}}$ is invertible, the first claim now follows.

To prove the second claim we first suppose that $R=\ZZ$. If $x,y\in
S(\lambda)$ then, by extending scalars, we may assume that $x$ and $y$
are elements of $S(\lambda)_{\Q(q)}=S(\lambda)_\ZZ\otimes\Q(q)$. Now
$S(\lambda)_{\Q(q)}\cong\pi\big(S(\lambda)_{\Q(q)}\big)$ is
irreducible so, up to a scalar, there is a unique associative bilinear form on
$S(\lambda)_{\Q(q)}$. To determine this scalar it is enough to compare
the two inner products
on $z_\lambda$ and $\pi(z_\lambda)$. Using associativity,
\begin{align*}
\<z_\lambda,z_\lambda\>_\lambda
     &=\<x_\lambda T_{w_\lambda} y_{\lambda'},
         x_\lambda T_{w_\lambda} y_{\lambda'}\>_\lambda
      =\<x_\lambda T_{w_\lambda} y_{\lambda'}^2,
           x_\lambda T_{w_\lambda}\>_\lambda\\
     &=q^{-\binom{k+1}2} [k+1]_q^!\<x_\lambda T_{w_\lambda} y_{\lambda'},
                         x_\lambda T_{w_\lambda}\>_\lambda\\
     &=q^{-\binom{k+1}2}[k+1]_q^!\sum_{v\in\Sym_{\lambda'}}(-q)^{-\len(v)}
        \<x_\lambda T_{w_\lambda v},x_\lambda T_{w_\lambda}\>_\lambda\\
     &=q^{\len(w_\lambda)-\binom{k+1}2}[k+1]_q^!\\
     &=q^{\frac k2(2n-3k-3)}[k+1]_q^!.
\intertext{Similarly,}
\big\<\pi(z_\lambda),\pi(z_\lambda)\big\>_{(k|n{-}k)}
     &=\<x_{(k|n{-}k)}y_{k+1}',x_{(k|n{-}k)}y_{k+1}'\>_{(k|n{-}k)}
      =q^{-k}[k+1]_q
\end{align*}
This proves that 
$\<x,y\>_\lambda=q^{\frac k2(2n-3k-1)}[k]_q^!\<\pi(x),\pi(y)\>_{(k|n{-}k)}$, for all 
$x,y\in S(\lambda)$ when $R=\ZZ$. The general case now follows by
specialization.
\end{proof}

\begin{Corollary}
Suppose that $\lambda=(n-k,1^k)$. Then $[k]_q^!$ divides
$\<x,y\>_\lambda$, for all $x,y\in S(\lambda)$.
\end{Corollary}

Let $S'(\lambda)=\pi\(S(\lambda)\)=x_{(k|n{-}k)}y_{k+1}'\H$. Then
$S'(\lambda)\cong S(\lambda)$ by the Proposition. We will work with
$S'(\lambda)$ in what follows rather than working
with $S(\lambda)$ directly. 

As a first step we need a basis of $S'(\lambda)$. For \textit{any}
$\lambda$-tableau $\t$ define
$$v_\t'=\pi(v_\t)=x_{(k|n{-}k)}y_{k+1}'T_{d(\t')}.$$

The Dipper--James basis of $S(\lambda)$, (\ref{DJ_basis}), combined
with Proposition~\ref{changing}, give us the following.

\begin{Corollary}
The module $S'(\lambda)$ is $R$--free with basis 
$\set{v_\t'|\t\in\Std(\lambda)}$.
\end{Corollary}

%
%
%

In order to exploit this basis we introduce another type of tableaux.
For our purposes we could get by using $(k,n-k)$--tableaux; however,
we use the notation from the theory of trivial source modules.

The {\sf diagram} of $(k|n{-}k)$ is the ordered pair of diagrams
$[k|n-k]=([k], [n-k])$.  A $(k|n{-}k)$--tableau is a bijection from
$[k|n-k]$ to $\{1,2,\dots,n\}$. Once again, we will think of a
$(k|n{-}k)$--tableau as being a labelling of $[k|n-k]$. Accordingly, we
will write a $(k|n{-}k)$--tableau as an ordered pair $(\a|\b)$, where
$\a$ and $\b$ are suitable labellings of the diagrams of the
partitions $(k)$ and $(n-k)$ respectively. We refer to $\a$ and $\b$
as the first and second components of $(\a|\b)$.

A $(k|n{-}k)$--tableau $(\a|\b)$ is {\sf (row) standard} if the
entries in $\a$ increase from left to right and the entries in $\b$
increase from left to right.  Let $\Std(k|n-k)$ be the set of standard
$(k|n{-}k)$--tableaux. For example, the standard $(1|3)$--tableaux are
$$\tab(1|234),\quad\tab(2|134),\quad\tab(3|124)\quad\And\quad
\tab(4|123).$$ Let $\t^{(k|n{-}k)}$ be the standard
$(k|n{-}k)$--tableau with $1,\dots,k$ entered in order, from left to
right, in the first omponent and the numbers $k+1,\dots,n$ in the
second. The first of the tableaux above is $\t^{(1|3)}$.

Two $(k|n-k)$--tableaux $(\a|\b)$ and $(\s|\t)$ are {\sf row
equivalent} if $\a$ and $\s$ contain the same entries up to reordering
(in which case, $\b$ and $\t$ also contain the same set of entries).
As with ordinary tableaux, the symmetric group acts from the right on
the set of $(k|n{-}k)$--tableaux. If $(\a|\b)$ is
a $(k|n{-}k)$--tableau we define $d(\a|\b)$ to be the unique
permutation such that $(\a|\b)=\t^{(k|n{-}k)}d(\a|\b)$. Then $(\a|\b)$
and $(\s|\t)$ are row equivalent if and only if $d(\a|\b)=wd(\s|\t)$
for some $w\in\Sym_{(k,n-k)}$. Consequently,
$$\D_{(k,n-k)}=\{\,d(\a|\b)\mid(\a|\b)\in\Std(k|n{-}k)\}.$$
So the standard $(k|n-k)$--tableaux index a basis of
$M(k|n{-}k)$. 

For furture reference, notice that if $(\a|\b)$ is a standard
$(k|n-k)$--tableau then $\len(d(\a|\b))$ is equal to the number of
pairs of integers $(i,j)$ where $i$ appears in $\a$, $j$ appears in
$\b$ and $i>j$. This follows because if $w\in\Sym_n$ then $\len(w)$ is
equal to the number of pairs $a<b$ with $i=a^w>j=b^w$, and the entries
in $\a$ are the images of $1,\dots,k$ under $d(\a|\b)$, whereas the
entries in $\b$ are the images of $k+1,\dots,n$.

For any $(k|n{-}k)$--tableau $(\a|\b)$ define
$x_{(\a|\b)}=x_{(k|n{-}k)}T_{d(\a|\b)}$. Here we do not assume that
$(\a|\b)$ is standard.  The following lemma is easily verified.

\begin{Lemma}\label{signed} Suppose that $0\le k<n$. 
\begin{enumerate}
\item $M(k|n{-}k)$ is free as an $R$--module with basis
$\{\,x_{(\a|\b)}\mid(\a|\b)\in\Std(k|n{-}k)\,\}$.
\item Suppose that $(\a|\b)\in\Std(k|n{-}k)$ and $1\le i<n$. Then
$$x_{(\a|\b)}T_i=\begin{cases}
  -x_{(\a|\b)},&\text{if $i$ and $i+1$ are both contained in $\a$,}\\
  qx_{(\a|\b)},&\text{if $i$ and $i+1$ are both contained in $\b$,}\\
  x_{(\a_i|\b_i)},&\text{if $i$ is in $\a$ and $i+1$ is in $\b$},\\
  qx_{(\a_i|\b_i)}+(q-1)x_{(\a|\b)},&\text{otherwise,}
\end{cases}$$
where $(\a_i|\b_i)=(\a|\b)r_i$.
\end{enumerate}
\end{Lemma}

The action of $\H$ on $M(k|n-k)$ is completely determined by (ii).

We now show how to write the basis $\{v_\t'\}$ of $S'(\lambda)$ in
terms of this basis of $M(k|n-k)$. To do this, if $\t$ is a
$\lambda$--tableau and $(\a|\b)$ is a $(k|n{-}k)$--tableau write
$(\a|\b)\prec\t$ if $(\a|\b)$ is standard and all of the entries in
$\a$ are contained in the first column of $\t$. Finally, if
$(\a|\b)\prec\t$ we set $I_\t(\a|\b)=i$, the index of $(\a|\b)$ in
$\t$, if the number in row $i$ of $\t$ does not appear in $\a$.

\begin{Lemma}\label{+lot}
Suppose that $\t$ is a standard $\lambda$--tableau. Then
$$v_\t'=\sum_{(\a|\b)\prec\t}
   (-1)^{k+1-I_\t(\a|\b)}q^{\len(d(\t'))-\len(d(\a|\b))}x_{(\a|\b)}.$$
\end{Lemma}

\begin{proof}
First consider $v_{\t_\lambda}$. Looking at the definitions we see
that
\begin{align*}
v_{\t_\lambda}'&=x_{(k|n-k)}y_{k+1}'
     =x_{\t^{(k|n-k)}}\Big(1-q^{-1}T_k+q^{-2}T_{k,k-1}
                 -\dots+(-q)^{-k}T_{k,1}\Big)\\
    &=\sum_{(\a|\b)\prec\t_\lambda}(-q)^{-\len(d(\a|\b))} x_{(\a|\b)}.
\end{align*}
As $\len(d(\t_\lambda)')=\len(d(\t^{\lambda'}))=0$ and
$\len(d(\a|\b))=k+1-I_{\t_\lambda}(\a|\b)$, when
$(\a|\b)\prec\t_\lambda$, the Lemma follows in this case.

Now suppose that $\t$ is an arbitrary standard $\lambda$--tableaux.
If $\t\ne\t_\lambda$ then we can find another standard
$\lambda$--tableau $\s$ and an integer $i$ in the first column of $\s$
such that~$\t=\s r_i$ and $\len(d(\t))=\len(d(\s))-1$. (That is,
$\t\gedom\s$ where $\gedom$ is the dominance order on tableaux; see,
for example,~\cite{M:IHA}.) Therefore, by induction,  
$$v_\t'=v_\s'T_i=\sum_{(\a|\b)\prec\s}
   (-1)^{k+1-I_\s(\a|\b)}q^{\len(d(\s'))-\len(d(\a|\b))}x_{(\a|\b)}T_i.$$
Since $\s$ and $\t$ are standard, $i$ is in the first column of $\s$
and the first row of $\t$ and~$i+1$ is in the first row of $\s$ and
the first column of $\t$.  Therefore, if $(\a|\b)\prec\s$ then the
entries in the first component of $(\a|\b)r_i$ are still in increasing
order and the entries in the second component are in increasing order
unless $i$ and~$i+1$ both appear in~$\b$. So,
$\len(d(\a|\b)r_i)=\len(d(\a|\b))+1$ and by Lemma \ref{signed}(ii) we
have
$$x_{(\a|\b)}T_i=\begin{cases}
      qx_{(\a|\b)},&\If i\And i+1\text{ both appear in $\b$},\\
       x_{(\a|\b)r_i},&\Otherwise.
\end{cases}$$
In the first case, when $i$ and $i+1$ both appear in $\b$, we have
that $(\a|\b)\prec\t$. Also, that
$\len(d(\t'))-\len(d(\a|\b))=\len(d(\s'))-\len(d(\a|\b))+1$ and
$I_\t(\a|\b)=I_\s(\a|\b)$ so $x_{(\a|\b)}$ has the required
coefficient in $v_\t'$.

In the second case, $i$ appears in $\a$ and $i+1$ appears in $\b$, so 
$(\a r_i|\b r_i)=(\a|\b)r_i\prec\t$, $I_\t(\a r_i|\b r_i)=I_\s(\a|\b)$
and $\len(d(\t'))-\len(d(\a|\b))=\len(d(\s'))-\len(d(\a|\b))+1$.
Hence, once again, $x_{(\a|\b)r_i}$ has the predicted coefficient in $v_\t'$.

As there are exactly $k$ standard $(k|n-k)$--tableaux $(\a|\b)$ 
satisfying $(\a|\b)\prec\t$, this completes the proof.
\end{proof}

In order to compute the elementary divisors of $S(\lambda)$ we need a
second basis of~$S'(\lambda)$. Let
$$x_{n-k}=1+T_1+\dots+T_{1,n-k-1}=\sum_{j=0}^{n-k-1}T_{1,j}.$$
(Note that $r_{1,0}=1$.) As with $y_{k+1}'$, we have
$x_{(n-k,1^k)}=x_{(1,n-k-1,1^{k-1})}x_{n-k}$. Now, for any standard
$(n-k,1^k)$--tableau $\t$ we define
$$w_\t'=\begin{cases}
          v_{\t(1,n)}',&\text{\ if $n$ appears in row 1 of $\t$},\\
          v_{\t^\lambda}'x_{n-k}T_{d(\t)},&\Otherwise.
\end{cases}$$
We remark that it is not obvious that the set of elements
$\set{w_\t'|\t\in\Std(\lambda)}$ is a basis of $S'(\lambda)$. We will
prove this below.

Lemma~\ref{+lot} gives an explicit description of the basis
$\{v_\t'\}$. We need to do the same for the basis $\{w_\t'\}$, and for
this we need some more notation. If $\t$ is a standard
$\lambda$--tableau let $\t^*=\t(1,n)$. If $(\a|\b)$ is a
$(k|n-k)$--tableau write $(\a|\b)\prec_n\t$ if $(\a|\b)\prec\t$ and
$n$ is contained in $\b$. Finally, if $n$ appears in the first row of
$\t$ then we define $(\a_\t^*|\b_\t^*)$ to be the unique standard
$(k|n-k)$--tableau such that $\(\a_\t^*|\b_\t^*)\prec\t^*$ and $n$
appears in $\b_\t^*$. So $\(\a_\t^*|\b_\t^*)\prec\t^*$ and
$\(\a_\t^*|\b_\t^*)\prec\t$.

\begin{Lemma}
Suppose that $\t$ is a standard $\lambda$--tableau and that $n$
appears in the first row of~$\t$.  Then
$$w_\t'=(-1)^kq^{2n-2k-3}x_{(\a_\t^*|\b_\t^*)}
        +\sum_{(\a|\b)\prec_n\t^*}r_{\a\b}x_{(\a|\b)},$$
for some scalars $r_{\a\b}\in\ZZ$. 
\label{+hot}
\end{Lemma}

\begin{proof}
We now argue by downwards induction on $\t$ beginning with
$\t=\t_\lambda$, this is an unpleasant calculation. Now, 
\begin{align*}
w_{\t_\lambda}'&=x_{(k|n{-}k)}y_{k+1}'T_{1,n-1}T_{n-2,1}
               =(-1)^kx_{(k|n{-}k)}y_{k+1}'T_{k+1,n-1}T_{n-2,1}\\
\intertext{since $x_{(k|n{-}k)}y_{k+1}'=x_{(1^k,n-k)}y_{(k+1,1^{n-k-1})}$.
Therefore, using the definitions together with the braid relations,}
w_{\t_\lambda}'
           &=(-1)^kx_{(k|n{-}k)}y_{k+1}'T_{n-1}\dots T_{k+2}T_{k+1,n-1}T_{k,1}\\
           &=(-1)^kx_{(k|n{-}k)}T_{n-1}\dots T_{k+2}y_{k+1}'T_{k+1,n-1}T_{k,1}\\
           &=(-1)^kq^{n-k-2}x_{(k|n{-}k)}y_{k+1}'T_{k+1,n-1}T_{k,1}\\
\end{align*}
\begin{align*}
           &=(-1)^kq^{n-k-2}x_{(k|n{-}k)}\Big\{1
             +\sum_{j=1}^k(-q)^{j-k-1}T_{k,j}\Big\}T_{k+1,n-1}T_{k,1}\\
           &=(-1)^kq^{n-k-2}x_{(k|n-k)}\Big\{q^{n-k-1}
             +\sum_{j=1}^k(-q)^{j-k-1}T_{k,j}T_{k+1,n-1}\Big\}T_{k,1}\\
           &=(-1)^kq^{n-2k-3}x_{(k|n-k)}\Big\{q^{n}T_{k,1}
             +\sum_{j=1}^k(-q)^{j}T_{k,j}T_{k+1,1}T_{k+2,n-1}\Big\}\\
           &=(-1)^kq^{n-2k-3}x_{(k|n-k)}\Big\{q^{n}T_{k,1}
           +\sum_{j=1}^k(-q)^{j}T_{k+1,1}T_{k+1,j+1}T_{k+2,n-1}\Big\}\\
           &=(-1)^kq^{n-2k-3}x_{(k|n-k)}\Big\{q^{n}T_{k,1}
           -\sum_{j=1}^k(-q)^{j+1}T_{k,1}T_{k+1,j+1}T_{k+2,n-1}\Big\}\\
           &=(-1)^kq^{n-2k-3}x_{(k|n-k)}T_{k,1}\Big\{q^{n}
             -\sum_{j=2}^{k+1}(-q)^{j}T_{k+1,j}T_{k+2,n-1}\Big\}.
\end{align*}
Now, 
$\t^{(k|n-k)}r_{k,1}=\big(\,\widetab(2&&\cdots&&k+1)\,\big|\,
                            \widetab(1&&k+2&&\cdots&&n)\,\big)
                    =(\a_{\t_\lambda}^*|\b_{\t_\lambda}^*)$ and,
consequently,
$\t^{(k|n-k)}r_{k,1}r_{k+1,j}r_{k+2,n-1}
    =\big(\,\widetab(2&&\cdots&&j-1&&j+1&&\cdots&&k+1&&n)\,\big|\,%
            \widetab(1&&j&&k+2&&\cdots&&n-1)\,\big),$
for $j=2,\dots,k+1$. This completes the proof for $w_{\t_\lambda}'$. 

Now suppose that $\t$ is an arbitrary standard $\lambda$--tableau
which has $n$ in its first row. Then $d(\t')\in\Sym_{(1,n-2,1)}$ so
$d(\t')$ and $(1,n)$ commute and
$\len(d(\t')(1,n))=\len(d(\t'))+\len(1,n)$.  Therefore,
$w_\t'=w_{t_\lambda}'T_{d(\t')}$. To complete the proof now argue by
induction, as in the proof of Lemma~\ref{+lot}; we leave the details
to the reader. (Indeed, this shows that $r_{\a\b}=\pm q^a$ for some
integer $a$.)
\end{proof}

For convenience we now write $\<\ ,\ \>=\<\ ,\ \>_{(k|n{-}k)}$. In
terms of the standard basis of $M(k|n-k)$, the bilinear form 
$\<\ ,\ \>$ on $M(k|n-k)$ is determined by
$$
\<x_{(\a|\b)}, x_{(\s|\t)}\> \; =\; 
\begin{cases}
q^{\len(d(\a|\b))},& \If (\a|\b)=(\s|\t),\\
0,                 & \Otherwise,
\end{cases}
$$
for standard $(k|n-k)$--tableaux $(\a|\b)$ and $(\s|\t)$.

\begin{Corollary}\label{top half}
Suppose that $\s$ and $\t$ are standard $\lambda$--tableaux which have
$n$ in their first row. Then
$$\<w_\s',v_\t'\>=\begin{cases}q^{2n-2k-3+\len(d(\t'))},&\If\s=\t,\\
                               0,&\Otherwise.
\end{cases}$$
\end{Corollary}

\begin{proof}By Lemma~\ref{+lot} and Lemma~\ref{+hot} we have
\begin{align*}
v_\s'&=\sum_{(\a_s|\b_s)\prec\s}
    (-1)^{k+1-I_\s(\a_s|\b_s)}q^{\len(d(\s'))-\len(d(\a_s|\b_s))}
              x_{(\a_s|\b_s)}\\
\intertext{and, by Lemma \ref{+hot},} 
w_\t'&=(-1)^kq^{2n-2k-3}x_{(\a_\t^*|\b_\t^*)}
        +\sum_{(\a_t|\b_t)\prec_n\t^*}r_{\a_t\b_t}x_{(\a_t|\b_t)}.
\end{align*}
Now, all of the tableaux appearing in $w_\t'$ have $1$
appearing in their second component. In contrast, the only tableau in
$v_\s'$ which has $1$ in its second component is the tableau
$(\a_\s^*|\b_\s^*)$. Therefore,
\begin{align*}
\<w_\s',v_\t'\>
    &=(-1)^{2k+1-I_\s(\a_s^*|\b_s^*)}
      q^{2n-2k-3+\len(d(\s'))-\len(d(\a_\s^*|\b_\s^*))}
             \<x_{(\a_\s^*|\b_\s^*)},x_{(\a_\t^*|\b_\t^*)}\>\\
    &=\begin{cases}
          (-1)^{1+I_\t(\a_\t^*|\b_\t^*)}
           q^{2n-2k-3+\len(d(\t'))},&\If\s=\t,\\  
           0,&\Otherwise.
      \end{cases}
\end{align*}
Finally, the sign vanishes when $\s=\t$ because $I_\t(\a_\t^*|\b_\t^*)=1$ .
\end{proof}

We need one more result before we can produce the elementary divisors
of $S(\lambda)$.

\begin{Lemma}Let $(\a_\lambda^+|\b_\lambda^+)$ be the unique standard
$(k|n-k)$--tableau which has the numbers $n-k+1,\dots,n$ in
$\a_\lambda^+$. Then
\begin{align*} w_{\t^\lambda}'
=(-1)^kq^{\len(w_{\lambda'})}\Big\{&
        q^{-\len(d(\a_\lambda^+|\b_\lambda^+))}[n-k]_q 
             x_{(\a_\lambda^+|\b_\lambda^+)}\\
  &\quad +\sum_{\substack{(\a|\b)\prec\t^\lambda\\
                         (\a|\b)\ne(\a_\lambda^+|\b_\lambda^+)}}
        (-1)^{1-I_{\t^\lambda}(\a|\b)}q^{-\len(d(\a|\b))}
        \sum_{j=0}^{n-k-1}x_{(\a|\b)r_{1,j}}\Big\}.
\end{align*}
\label{w_t}\end{Lemma}

\begin{proof}
By definition $w_{\t^\lambda}'=v_{\t^\lambda}'x_{n-k}$. Also,
$d((\t^\lambda)')=d(\t_{\lambda'})=w_{\lambda'}$ so, by 
Lemma~\ref{+lot},
\begin{align*}w_{\t^\lambda}'
 &=v_{\t^\lambda}'x_{n-k}=\sum_{(\a|\b)\prec\t^\lambda}
               (-1)^{k+1-I_{\t^\lambda}(\a|\b)}
               q^{\len(w_{\lambda'})-\len(d(\a|\b))}x_{(\a|\b)}x_{n-k}\\
 &=\sum_{(\a|\b)\prec\t^\lambda}
    (-1)^{k+1-I_{\t^\lambda}(\a|\b)}q^{\len(w_{\lambda'})-\len(d(\a|\b))}
        x_{(\a|\b)}\Big(1+T_1+\dots+T_{1,n-k-1}\Big).
\end{align*}
Let $(\a|\b)$ be one of the tableaux appearing in this sum. If
$(\a|\b)\ne(\a_\lambda^+|\b_\lambda^+)$ then $1$ is contained in $\a$
and all of the numbers $2,3,\dots,n-k$ are contained in $\b$.
Therefore, $(\a|\b)r_{1,j}$ is standard and
$x_{(\a|\b)}T_{1,j}=x_{(\a|\b)r_{1,j}}$, for $0\le j\le n-k-1$. On the
other hand,
$x_{(\a_\lambda^+|\b_\lambda^+)}T_{1,j}=q^jx_{(\a_\lambda^+|\b_\lambda^+)}$,
for $0\le j\le n-k-1$.  This completes the proof of the Lemma.
\end{proof}

This result has two useful Corollaries.

\begin{Corollary}\label{zeros}
Suppose that $\t\ne\t^\lambda$ is a standard $(k|n-k)$--tableau. Then
$$\<w_{\t^\lambda}',v_\t'\>=0.$$
\end{Corollary}

\begin{proof}
By Lemma~\ref{w_t}, if $x_{(\a|\b)}$ appears in $w_{\t^\lambda}'$ then
all but one of the entries in $\a$ are contained in
$\{1,n-k+1,\dots,n\}$.  On the other hand, by Lemma~\ref{+lot}, if
$x_{(\a|\b)}$ appears in $v_\t'$ then all of the entries in $\a$ are
contained in the first column of~$\t$. 

Suppose now that $\t\ne\t^\lambda$. Then, by the last paragraph,
$x_{(\a_\lambda^+|\b_\lambda^+)}$ cannot appear in~$v_\t'$ and the
only way that the inner product $\<w_{\t^\lambda}',v_\t'\>$ can be
non--zero is if the set of numbers in the first column of $\t$ is of
the form $T=\{1,j,n-k+1,\dots,n\}\setminus\{m\}$, for some integers
$j$ and $m$ with $1<j\le n-k$ and $n-k<m\le n$. Let $(\a|\b)$ be the
standard $(k|n-k)$--tableau whose first component contains exactly the
numbers in~$T\setminus\{j\}$ and let $(\s|\t)=(\a|\b)r_{1,j-1}$. Then
$(\a|\b)\prec\t^\lambda$, $I_\t(\a|\b)=2$ and $I_\t(\s|\t)=1$. Also
$\len(d(\s|\t))=\len(d(\a|\b))+j-1$, so
$x_{(\a|\b)}T_{1,j-1}=x_{(\s|\t)}$. Therefore, by Lemma~\ref{+lot} and
Lemma~\ref{w_t} and the remarks above,
\begin{align*}
v_\t'&=(-1)^kq^{\len(d(\t'))}
    \Big(q^{-\len(d(\s|\t))}x_{(\s|\t)}-q^{-\len(d(\a|\b))}x_{(\a|\b)}\Big)
           +\text{\ other standard terms}\\
\intertext{and}
w_{\t^\lambda}'&=q^{\len(w_{\lambda'})-\len(d(\a|\b))}
                     \Big(x_{(\s|\t)}+x_{(\a|\b)}\Big)
+\text{other standard terms,}
\end{align*}
where none of the ``other standard terms'' appear both in $v_\t'$ and
in $w_{\t^\lambda}'$. Consequently, $\<w_{\t^\lambda}',v_\t'\>=0$.
Hence, $\<w_{\t^\lambda}',v_\t'\>=0$ whenever $\t\ne\t^\lambda$ as
claimed.
\end{proof}

\begin{Corollary}\label{n-diagonals}
Suppose that $\t$ is a standard $(n-k,1^k)$--tableau and that
$n$ does not appear in the first row of $\t$. Then
$\<w_\t',v_\t'\>=q^{k(n-k-2)}[n]_q.$
\end{Corollary}

\begin{proof}
Recall that if $\t$ is a standard $\lambda$--tableau
then $d(\t')d(\t)^{-1}=w_{\lambda'}$, with the lengths adding; this is
well--known and is easily proved by induction on the dominance order
for tableaux. Therefore, 
\begin{align*}
\<w_\t',v_\t'\>&=\<w_{\t^{\lambda}}'T_{d(\t)},x_{(k|n{-}k)}y_{k+1}'T_{d(\t')}\>
          =\<w_{\t^{\lambda}}',x_{(k|n{-}k)}y_{k+1}'T_{d(\t')}T_{d(\t)}^*\>\\
         &=\<w_{\t^{\lambda}}',x_{(k|n{-}k)}y_{k+1}'T_{d({\t^{\lambda}}')}\>
          =\<w_{\t^\lambda}',v_{\t^\lambda}'\>.
\end{align*}
Hence, it is enough to consider the case where $\t=\t^\lambda$.

Suppose that $\t=\t^\lambda$. Then, by Lemma~\ref{+lot} and
Lemma~\ref{w_t},
\begin{align*}
\<w_{\t^\lambda}',v_{\t^\lambda}'\>
     &=q^{2\len(w_{\lambda'})}\Big\{
                        q^{-\len(d(\a_\lambda^+|\b_\lambda^+))}[n-k]_q
          +\sum_{\substack{(\a|\b)\prec\t^\lambda\\
                           (\a|\b)\ne(\a_\lambda^+|\b_\lambda^+)}}
                        q^{-\len(d(\a|\b))}\Big\}.
\intertext{Using the remarks before Lemma~\ref{signed} it is not hard
to see that $\len(d(\a_\lambda^+|\b_\lambda^+))=k(n-k)$ and that
$\len(d(\a|\b))=(k-1)(n-k)+2-I_\t(\a|\b)$, whenever
$(\a|\b)\prec\t^\lambda$ and $(\a|\b)\ne(\a_\lambda^+|\b_\lambda^+)$. 
Therefore,}
\<w_{\t^\lambda}',v_{\t^\lambda}'\>
     &=q^{2\len(w_{\lambda'})}\Big\{
         q^{-k(n-k)}[n-k]_q+\sum_{i=2}^{k+1}q^{-(k-1)(n-k)-2+i}\Big\}\\
     &=q^{2\len(w_{\lambda'})-k(n-k)}\Big\{
                       [n-k]_q+q^{n-k}\sum_{j=0}^{k-1}q^j\Big\}\\
     &=q^{2\len(w_{\lambda'})-k(n-k)}[n]_q.
\end{align*}
As $\len(w_{\lambda'})=k(n-k-1)$ the result follows.
\end{proof}

Finally, we can prove the main result of this section.

\begin{Proposition}\label{eds}
Suppose that $\lambda = (n-k,1^k)$, for some $k$ with $0\le k<n$.
Then the Gram matrix $G(\lambda)$ of $S(\lambda)$ is divisibly
diagonalizable over $\ZZ$ with $\binom{n-2}{k}$ elementary divisors
equal to $[k]_q^!$ and with the remaining $\binom{n-2}{k-1}$
elementary divisors being equal to $[k]_q^![n]_q$.
\end{Proposition}

\begin{proof}By Proposition~\ref{changing} the Gram matrix
$G(\lambda)$ of $S(\lambda)$ is equal to $[k]_q^!$ times the Gram
matrix of $S'(\lambda)$. Therefore, by Lemma~\ref{diag} it is enough
to show that there is an invertible diagonal matrix $D$ such that
$$G'(\lambda)=\Big(\<w_\s',v_\t'\>\Big)_{\s,\t\in\Std(k|n-k)}
     = D\cdot\left(\begin{array}{*4c} 
                     I & *\\ 
                     0 & [n]_qU
           \end{array}\right),$$ 
where $I$ is a $\binom{n-2}k\times\binom{n-2}k$ identity matrix and
$U$ is a $\binom{n-2}{k-1}\times\binom{n-2}{k-1}$ upper triangular matrix
with $1$'s down its diagonal. Here we order the rows and columns
lexicographically with respect to the entries in the first column of
$\s$ and $\t$. Because $D$ is invertible its non--zero entries are all
of the form $\pm q^m$, for some integer $m$.

By Corollary~\ref{top half}, the rows of $G'(\lambda)$ which are
indexed by those tableaux which have $n$ in their first row have the
required form. This accounts for the identity matrix in the top half
of the Gram matrix $G'(\lambda)$.

Next, suppose that $\s$ is a standard
$(k|n-k)$--tableau and that $n$ does not appear in the first row of
$\s$. If $\s =\t^\lambda$ then $\<w_\s',v_\t'\>=0$, for all $\t\ne\s$,
by Corollary~\ref{zeros}.
If $\s\ne\t^\lambda$ then there exists an integer $i$, $1\le i<n$, such
that $\len(d(\s)r_i)<\len(d(\s))$. Therefore,
$$
\<w_\s',v_\t'\>=\<w_{\s r_i}'T_i,v_\t'\>
                 =\<w_{\s r_i}',v_\t'T_i\>.
$$
It follows that $\<w_\s',v_\t'\>=0$ if $\t$ appears before $\s$ in our
chosen ordering of $\Std(\lambda)$. Finally, if $\t$ does not appear
before $\s$ then $[n]_q$ divides $\<w_\s',v_\t'\>$ by
Corollary~\ref{n-diagonals}.
\end{proof}

Notice, in particular, that the Gram matrix calculation in the proof of
the Proposition implies that $\set{w_\t'|\t\in\Std(\lambda)}$ is
indeed a basis of $S'(\lambda)$. 

Proposition~\ref{eds} also gives the decomposition numbers of
$S(\lambda)$ (by inducing the corresponding principal indecomposable
modules); however, as these are already known we leave these as an
exercise for the reader. We will, however, give one application of
this result.

Let $\pi\map{S(\lambda)}S'(\lambda)$ be the isomorphism of
Proposition~\ref{changing} and for each standard $\lambda$--tableau
$\t$ let $w_\t=\pi^{-1}(w_\t')$. Then $\set{w_\t|\t\in\Std(\lambda)}$
is a basis of $S(\lambda)$.  Then, in the case where $S(\lambda)$ is
not irreducible, the proof of Proposition~\ref{eds} also gives a basis
for the simple module $D(\lambda)$. More precisely, we have the
following.

\begin{Corollary}
Suppose that $R$ is a field, that $[k]_q^!\neq 0$ and that $[n]_q = 0$. Then
$S(\lambda)$ is not irreducible and a basis of
$D(\lambda) = S(\lambda)/\big(S(\lambda)^\perp\cap S(\lambda)\big)$ is given by 
$$
\set{w_\t + \big(S(\lambda)^\perp\cap S(\lambda)\big)|%
               \t\in\Std(\lambda)\And\text{$n$ in first row of $\t$}}\;,
$$
and a basis of $S(\lambda)^\perp\cap S(\lambda)$ is given by
$
\set{w_\t |\t\in\Std(\lambda)\And\text{$n$ is in first row of $\t$}}.
$
\end{Corollary}

\section{Some counterexamples}
\label{SecCounterEx}

Let $R = \Z[q,q^{-1}]$. We write the $m^\th$ cyclotomic polynomial in
$q$ as $\Phi_m = \Phi_m(q)$.

Andersen remarked that in general the Gram matrix $G(\lambda)$ is not
diagonalizable \cite[Remark~5.11]{Andersen:Tilting}. We give two
examples of this kind.

Note that $G(\lambda)$ is divisibly diagonalizable over
$\Z_{(p)}[q,q^{-1}]$ for all but finitely many primes $p$. In fact, it
suffices to exclude the primes occurring in the denominators of the
entries of the matrices used to diagonalize $G(\lambda)$ over
$\Q[q,q^{-1}]$.

We record the elementary divisors in ``jump notation''. That is, we write
\[
\lraa{f_1} m_1 \lraa{f_2} m_2 \lraa{f_3} m_3 \lraa{f_4} \cdots \lraa{f_s} m_s
\]
to indicate that the matrix has the elementary divisor $f_1$ with
multiplicity $m_1$, the elementary divisor $f_1 f_2$ with multiplicity
$m_2$, \dots, and the elementary divisor $f_1\cdots f_s$ with
multiplicity $m_s$. 

\begin{Example}
\rm
Let $\lambda = (3,3,2)$. The elementary divisors of $G(3,3,2)$ over $\Q[q,q^{-1}]$ are given by
\[
\lraa{\Phi_2^2} 1 \lraa{\Phi_4} 20\lraa{\Phi_3\Phi_5} 20 \lraa{\Phi_4} 1\; ;
\]
over $\F_2[q,q^{-1}]$ they are given by
\[
\lraa{\Phi_2^3} 1 \lraa{\Phi_2} 20\lraa{\Phi_3\Phi_5} 20 \lraa{\Phi_2} 1\; ;
\]
and, putting $q = 1$, over $\Z$ they are given by
\[
\lraa{2^3} 21\lraa{3\cdot 5} 21 \; .
\]
We claim that $G(3,3,2)$ is not diagonalizable over
$\Z_{(2)}[q,q^{-1}]$. To see this suppose that it is diagonalizable.
Then, considered as an element of $\Z_{(2)}[q,q^{-1}]$, any resulting
diagonal entry must contain the factor $(q+1)$ with exponent $2$.
Considered as an element of $\F_2[q,q^{-1}]$ the factor $(q+1)$ can
occur only with even exponent in such a diagonal entry. But this is
not the case, so we have a contradiction.

This claim in particular implies that $G(3,3,2)$ is not diagonalizable
over $\Z[q,q^{-1}]$.

We remark that the comparison of the elementary divisors over $\Q[q,q^{-1}]$ and over 
$\Z$ yields a contradiction to diagonalizability over $\Z_{(2)}[q,q^{-1}]$, too.
\end{Example}

\begin{Example}
\rm
Let $\lambda = (4,2,1,1)$. The elementary divisors of $G(4,2,1,1)$ over $\Q[q,q^{-1}]$ are given by
\[
\lraa{\Phi_2} 14\lraa{\Phi_2} 1 \lraa{\Phi_4} 30\lraa{\Phi_7} 30\lraa{\Phi_4} 1 \lraa{\Phi_2} 14\; ;
\]
over $\F_2[q,q^{-1}]$ they are given by
\[
\lraa{\Phi_2} 14\lraa{\Phi_2^2} 1 \lraa{\Phi_2} 30\lraa{\Phi_7} 30\lraa{\Phi_2} 1 \lraa{\Phi_2^2} 14\; ;
\]
over $\F_3[q,q^{-1}]$ they are given by
\[
\lraa{\Phi_2} 13\lraa{\Phi_2} 2 \lraa{\Phi_4} 30\lraa{\Phi_7} 30\lraa{\Phi_4} 2 \lraa{\Phi_2} 13\; ;
\]
and, putting $q = 1$, over $\Z$ they are given by
\[
\lraa{2} 14\lraa{2^2} 31 \lraa{7} 31\lraa{2^2} 14\; .
\]
We claim that $G(4,2,1,1)$ is not diagonalizable over
$\Z_{(2)}[q,q^{-1}]$. Again, by way of contradiction suppose that it
is diagonalizable. In $\F_2[q,q^{-1}]$, $14$ of the resulting diagonal
entries contain the factor $(q+1)$ with exponent $1$. Therefore, in
$\Z_{(2)}[q,q^{-1}]$, $14$ of them contain the factor $(q+1)$ with
exponent $1$ and the factor $(q^2 + 1)$ with exponent $0$. Similarly,
in $\F_2[q,q^{-1}]$, $14$ of the resulting diagonal entries contain
the factor $(q+1)$ with exponent $7$. Thus in $\Z_{(2)}[q,q^{-1}]$,
$14$ of them contain the factor $(q+1)$ with exponent $3$ and the
factor $(q^2 + 1)$ with exponent $2$. Hence in $\F_2[q,q^{-1}]$, no
other diagonal entry can contain $(q+1)$ with odd exponent. But in
$\F_2[q,q^{-1}]$, there is a diagonal entry containing $(q+1)$ to the
power $3$ and another containing it to the power $5$ so, again, we have a
contradiction.

We claim that $G(4,2,1,1)$ is not diagonalizable over $\Z_{(3)}[q,q^{-1}]$. Assume it to be diagonalizable. In $\Z_{(3)}[q,q^{-1}]$, $14$ of the resulting
diagonal entries contain $(q+1)$ with exponent $1$. This contradicts the fact that in $\F_3[q,q^{-1}]$, only $13$ of them contain $(q+1)$ with 
exponent $1$.

Both claims independently imply that $G(4,2,1,1)$ is not diagonalizable over $\Z[q,q^{-1}]$.

We remark that the comparison of the elementary divisors over $\Q[q,q^{-1}]$ and over 
$\Z$ yields a contradiction to diagonalizability over $\Z_{(2)}[q,q^{-1}]$, too.
\end{Example}

Finally, we give a (non-exhaustive) list of elementary divisors of
some divisibly diagonalizable Gram matrices for non-hooks, calculated
using {\sc Gap 3} \cite{Gap} and {\sc Magma} \cite{Magma}. We omit the
respective conjugate partition; compare (\ref{CorDiagbar}).

\begin{footnotesize}
$$\begin{array}{lll}\toprule
n& \lambda&\text{Elementary divisors of\ } G(\lambda) \\\midrule
4&(2,2)   &\lraa{\Phi_2} 1\lraa{\Phi_3} 1 \\[4pt]\hline
5&(3,2)   &\lraa{1} 1\lraa{\Phi_3} 3\lraa{\Phi_4} 1 \\[4pt]\hline
6&(4,2)   &\lraa{1} 4\lraa{\Phi_4} 1 \lraa{\Phi_2} 3\lraa{\Phi_5} 1 \\ 
 &(3,3)   &\lraa{\Phi_2} 1\lraa{\Phi_3} 3\lraa{\Phi_4} 1 \\
 &(3,2,1) &\lraa{1} 4\lraa{\Phi_3} 4\lraa{\Phi_5} 4\lraa{\Phi_3} 4 \\[4pt]\hline
7&(5,2)   &\lraa{1} 8\lraa{\Phi_5} 5\lraa{\Phi_3\Phi_6} 1 \\
 &(4,3)   &\lraa{1} 1\lraa{\Phi_3} 7\lraa{\Phi_4} 5\lraa{\Phi_5} 1 \\
 &(3,3,1) &\lraa{\Phi_2} 6 \lraa{\Phi_3} 2\lraa{\Phi_5} 12\lraa{\Phi_4} 1\\[4pt]\hline
8&(6,2)   &\lraa{1} 13\lraa{\Phi_3\Phi_6} 1 \lraa{\Phi_2} 5\lraa{\Phi_7} 1 \\
 &(5,3)   &\lraa{1} 8\lraa{\Phi_4} 6\lraa{\Phi_2} 7\lraa{\Phi_5} 6\lraa{\Phi_6} 1 \\
 &(4,4)   &\lraa{\Phi_2} 1\lraa{\Phi_3} 7\lraa{\Phi_4} 5\lraa{\Phi_5} 1 \\[4pt]\hline
9&(7,2)   &\lraa{1} 19 \lraa{\Phi_7} 7 \lraa{\Phi_4\Phi_8} 1 \\
 &(6,3)   &\lraa{1}21\lraa{\Phi_5}19\lraa{\Phi_6} 1 \lraa{\Phi_3} 6\lraa{\Phi_7} 1\\
 &(5,4)   &\lraa{1} 1\lraa{\Phi_3} 15\lraa{\Phi_4} 18\lraa{\Phi_5} 7\lraa{\Phi_6} 1 
\\[4pt]\bottomrule
\end{array}$$
\end{footnotesize}

We do not know an example of a Gram matrix $G(\lambda)$ that is diagonalizable over
$\Z[q,q^{-1}]$, but not divisibly diagonalizable.

For a general partition $\lambda$, we can not decide whether
$G(\lambda)$ is diagonalizable over $\Z[q,q^{-1}]$.

\section*{Acknowledgements}
The authors would like to thank Adrian Williams for pointing out some
errors in an earlier version of this manuscript.

\let\em\it

\end{document}